\documentclass[microtype]{gtpart}

\usepackage{pdfsync}
\usepackage{comment}

\input xy
\xyoption{all}
\usepackage{amsmath}
\usepackage{amsthm}
\usepackage{amsfonts}
\usepackage{amscd}

\usepackage{epsfig}
\usepackage{enumerate}

\usepackage{amssymb}
\usepackage{marvosym}

\usepackage{hyperref}

\title{Homological Stability for automorphism groups of Raags}

\author{Giovanni Gandini} 
\givenname{Giovanni}
\surname{Gandini}
\address{}
\email{gandinigio@googlemail.com}
\urladdr{}

\author{Nathalie Wahl} 
\givenname{Nathalie}
\surname{Wahl}
\address{Department of Mathematical Sciences\\ University of
  Copenhagen\\\newline
Universitetsparken 5\\ DK-2100 Copenhagen\\ Denmark}
\email{wahl@math.ku.dk}
\urladdr{http://www.math.ku.dk/~wahl/}


\arxivreference{1510.06723}  
\arxivpassword{v5xv9}   

\renewcommand{\epsilon}{\varepsilon}

\newtheorem{thm}{Theorem}[section]
\newtheorem{proposition}[thm]{Proposition}
\newtheorem{prop}[thm]{Proposition}

\newtheorem{cor}[thm]{Corollary}
\newtheorem{lemma}[thm]{Lemma}
\newtheorem{lem}[thm]{Lemma}

\theoremstyle{definition}

\newtheorem{Def}[thm]{Definition}

\newtheorem{Th}{Theorem}

\theoremstyle{remark}
\newtheorem{remark}[thm]{Remark}

\renewcommand{\C}{\mathcal{C}}
\newcommand{\F}{\mathcal{F}}
\newcommand{\G}{\mathcal{G}}

\newcommand{\OO}{\mathcal{O}}

\newcommand{\al}{\alpha}

\newcommand{\Ga}{\Gamma}

\newcommand{\De}{\Delta}
\newcommand{\eps}{\epsilon}

\newcommand{\s}{\sigma}
\newcommand{\Si}{\Sigma}

\newcommand{\rar}{\longrightarrow}
\newcommand{\inc}{\hookrightarrow}
\newcommand{\sta}{\stackrel}
\newcommand{\arsim}{\sta{\simeq}{\rar}}
\newcommand{\minus}{\backslash}

\newcommand{\x}{\times}

\newcommand{\op}{\oplus}

\newcommand{\lgl}{\langle}
\newcommand{\rgl}{\rangle}

\newcommand{\GL}{\operatorname{GL}}

\newcommand{\del}{\partial}

\newcommand{\aut}{\operatorname{Aut}}

\newcommand{\link}{\operatorname{Link}}
\newcommand{\st}{\operatorname{Star}}
\newcommand{\Fix}{\operatorname{Fix}}
\newcommand{\Hom}{\operatorname{Hom}}

\begin{document}

\begin{abstract} We show that the homology of the automorphism group of a right-angled Artin group stabilizes under taking products with any right-angled Artin group.
\end{abstract}

\maketitle

\section*{Introduction}

It has been conjectured that, for any (finitely generated) discrete group $G$, the homology groups $H_i(\aut(G^{*n});\Z)$ and $H_i(\aut(G^n);\Z)$ should be independent of $n$, for $n\gg i$, generalizing the classical stability results for $\GL_n(\Z)$ and $\aut(F_n)$ when $G = \Z$. (See the conjectures \cite[Conjecture 1.4]{Hatcher-Wahl},\cite[Conjecture 5.16]{stab}, and the classical results in \cite{gln,hatcher,HatVog,Maazen,vdk}.)

The stabilization of  $H_i(\aut(G^{*n});\Z)$ for $i$ large has been shown to hold for most
groups by the main theorem of \cite{djament} and \cite[Corollary 1.3]{Hatcher-Wahl}.\footnote{\cite{djament} gives stability for
  $\aut(G^{*n})$ with $G$ any group with a finite free product decomposition (eg.~a finitely generated group) without $\Z$ factor, while \cite{Hatcher-Wahl} treats the cases with $G$ arising as fundamental groups of certain 3--manifolds, allowing $\Z$ factors in the free product decomposition.} The stabilization of 
$H_i(\aut(G^n);\Z)$ in contrast has so far only been known in two extreme cases: when $G$ is abelian and when $G$ has trivial center and does not factorize as a direct product. 
Indeed, in the first case $\aut(G^n)$ is isomorphic to  $\GL_n(\operatorname{End}(G))$, which is known to stabilize (see Proposition~\ref{abeliangp}), while in the second case, 
 the group $\aut(G^n)$ is isomorphic to $\aut(G)\wr \Si_n$ \cite{john}, a group that is also known to stabilize \cite[Proposition 1.6]{Hatcher-Wahl}.\footnote{Slightly more generally, for the second case, one can get stability for $\aut(G^n)$ for $G$ a product of certain such center-free groups using \cite{john}.} 
In the present paper, we verify that the second conjecture holds for $G$ any right-angled Artin group, possibly factorizable, possibly with a non-trivial center. This proves a first ``mixed case''  of the conjecture, which interpolates between the two previously known cases. 

\smallskip

A {\em right-angled Artin group} (or {\em RAAG}) is a group with a finite set of generators $s_1,\dots,s_n$ and relations that are commutation relations between the generators, i.e.~relations of the form $s_is_j=s_js_i$ for certain $i$'s and $j$'s. The extreme examples of RAAGs are the free groups $F_n$ if no commutation relation holds, and the free abelian groups $\Z^n$ if all commutation relations hold. Given any two RAAGs $A$ and $B$, their product is again a RAAG. 
We consider in the present paper the sequence of groups $G_n=\aut(A\x B^n)$ associated to $A$ and $B$, and the sequence of maps $$\s_n\co G_n=\aut(A\x B^n)\ \rar\  G_{n+1}=\aut(A\x B^{n+1})$$ 
taking an automorphism $f$ of $A\x B^n$ to the automorphism $f\x B$ of $A\x B^{n+1}$ leaving the last $B$ factor fixed. 
Note that when $A$ is the trivial group, the group $G_n=\aut(B^n)$ is a group as in the second conjecture above.
  
Our main result is the following:

\begin{Th}[Stability with constant coefficients]\label{maincst}
Let $A,B$ be any RAAGs. The map 
$$H_i(\aut(A\x B^n);\Z)\rar H_i(\aut(A\x B^{n+1});\Z)$$ 
induced by $\s_n$  is surjective  for all $i\le \frac{n-1}{2}$ and an isomorphism for $i\le \frac{n-2}{2}$. If $B$ has no $\Z$--factors, then surjectivity holds for 
$i\le \frac{n}{2}$ and injectivity for $i\le \frac{n-1}{2}$. 
\end{Th}

We prove this stability theorem using the general method developed by Randal-Williams and the second author in \cite{stab}. This method provides a more general stability result, namely stability in homology not only with constant coefficients $\Z$ as above, but also with both {\em polynomial} and {\em abelian} coefficients, and we establish our main result also in this level of generality as Theorem~\ref{allbounds}.  The following theorems are further special cases of Theorem~\ref{allbounds}:

\smallskip

Stability for $\aut(A\x B^n)$ with the (abelian) coefficients $H_1(\aut(A\x B^n))$ implies the following:

\begin{Th}[Stability for commutator subgroups]\label{mainab}
Let $A,B$ be any RAAGs and let $\aut'(A\x B^n)$ denote the commutator subgroup of $\aut(A\x B^n)$. The map 
$$H_i(\aut'(A\x B^n);\Z)\rar H_i(\aut'(A\x B^{n+1});\Z)$$ 
induced by $\s_n$ is surjective  for all $i\le \frac{n-2}{3}$ and an isomorphism for $i\le \frac{n-4}{3}$. If $B$ has no $\Z$--factors, then surjectivity holds for 
$i\le \frac{n-1}{3}$ and injectivity for $i\le \frac{n-3}{3}$. 
\end{Th}

An example of a polynomial coefficient system for the groups $\aut(A\x B^n)$ is the sequence of ``standard'' representations $H_1(A\x B^n)$, and stability with polynomial coefficients yields the following in that case: 

\begin{Th}[Stability with coefficients in the standard representation]\label{mainpol}
Let $A,B$ be any RAAGs. Then the map 
$$H_i(\aut(A\x B^n);H_1(A\x B^n))\rar H_i(\aut(A\x B^{n+1});H_1(A\x B^{n+1}))$$ 
 is surjective  for all $i\le \frac{n-2}{2}$ and an isomorphism for $i\le \frac{n-3}{2}$. If $B$ has no $\Z$--factors, then surjectivity holds for 
$i\le \frac{n-1}{2}$ and injectivity for $i\le \frac{n-2}{2}$. 
\end{Th}

To prove the above theorems, we show that right-angled Artin groups under direct product fit in the set-up of homogeneous categories developed in \cite{stab}, and recalled here in Section~\ref{sec1}. The main ingredient of stability is the high connectivity of certain semisimplicial sets $W_n(A,B)$ associated to the sequence of groups $\aut(A\x B^n)$. We define and study those semisimplicial sets in Section~\ref{complexes}, together with three closely related simplicial complexes $I_n(A,B)$, $SI_n(A,B)$ and $S_n(A,B)$. Sections~\ref{sec1} and \ref{complexes} are written in the general context of families of groups closed under direct product. In Section~\ref{RAAGsec}, we show that right-angled Artin groups admit a ``prime decomposition'' with respect to direct product, and we give a description of the automorphism group of such a group in terms of this decomposition. 
Section~\ref{connsec} then uses these results that are specific to RAAGs together with the complexes defined in Section~\ref{complexes} to prove that the semisimplicial sets $W_n(A,B)$ are highly connected. 
For the connectivity results, we  use {\em join complex} methods from
\cite{Hatcher-Wahl}, as well as an argument of  Maazen \cite{Maazen} for the case $B=\Z$.
Finally Section~\ref{stabsec} states the general stability result, which, given the connectivity result, is a direct application of the main result in \cite{stab}.

\subsection*{Acknowledgement} The second author was supported by the Danish National Sciences Research Council (DNSRC) and the European Research Council (ERC), as well as by the Danish National Research Foundation through the Centre for Symmetry and Deformation (DNRF92).

\section{Families of groups}\label{sec1}

We consider here families of  groups $\F$ which are
closed under direct product. 
We say that $\F$ satisfies {\em cancellation} if for all $A,B,C$ in $\F$, we have that 
$$A\x C\cong B\x C \ \ \Longrightarrow\ \ A\cong B.$$
Cancellation is not satisfied for the family of all finitely generated groups, see eg.~\cite[Section 3]{Hirshon1} or \cite{Hirshon0} for an example where cancellation with $\Z$ fails. 
Cancellation though holds for the family of all finitely generated abelian groups by their classification, the family of all finite groups \cite{Hirshon0}, or for the family of all right angled Artin groups as we will show in Section~\ref{RAAGsec}. 

\smallskip

Given a family of groups $\F$, we let $\G_\F$ denote its associated groupoid, namely the groupoid with objects the elements of $\F$ and morphisms 
all group isomorphisms. Let $e$ denote the trivial group. When $\F$ is closed under direct product, we have that 
$(\G_\F,\x,e)$ is a symmetric monoidal groupoid.

Recall from \cite[Section 1.1]{stab}, \cite[page 219]{Gra76} the category
$U\G_\F=\lgl \G_\F,\G_\F\rgl$ associated to $(\G_\F,\x,e)$: it has the
same objects as $\G_\F$,  namely the elements of $\F$, and morphisms
from $A$ to $B$ given as pairs $(X,f)$, where $X\in \F$ and $f\co X\x A\arsim B$ is an isomorphism, up to the equivalence relation that $(X,f)\sim (X',f')$ if there exists an isomorphism $\phi:X\to X'$ such that $f=f'\circ (\phi\x A)$.

Recall from \cite[Definition 1.2]{stab} that a monoidal category $(\C,\op,0)$ is called {\em homogeneous} if 0 is initial in $\C$ and for every $A,B$ in $\C$, the following two properties hold: 
\begin{itemize}
\item[{\bf H1}] $\Hom(A,B)$ is a transitive $\aut(B)$-set under post-composition;
\item[{\bf H2}] The map $\aut(A)\to \aut(A\oplus B)$ taking $f$ to $f\op B$ is injective with image $\Fix(B, A \op B)$,
\end{itemize}
where $\Fix(B, A \op B)$ is the set of $\phi\in\aut(A\op B)$ satisfying that $\phi\circ (\iota_A\op B)=\iota_A\op B$ in $\Hom(B,A\op B)$, for $\iota_A\co 0\to A$ the unique morphism. 

\begin{prop}\label{homogeneous} If $\F$ satisfies cancellation, then 
the category $U\G_{\F}$ is a symmetric monoidal homogeneous category whose underlying groupoid is $\G_\F$. 
\end{prop}

\begin{proof}
As $(\G_\F,\x,e)$ is symmetric monoidal, $U\G_\F$ is symmetric monoidal by \cite[Proposition 1.6]{stab}, and $e$ is initial in $U\G_\F$. 
We have that $\G_\F$ satisfies cancellation by assumption, and 
for any $A,B\in \F$, the map $\aut_{\G_\F}(A)\to \aut_{\G_\F}(A\x B)$ taking $f$ to $f\x B$ is injective. Then \cite[Theorem 1.8]{stab} implies that $U\G_\F$ is a homogeneous category. 
Finally, if $A\x B\cong e$, we must have $A=B=e$ and the unit $e$ has no non-trivial automorphisms. 
Hence $\G_\F$ satisfies the hypothesis of Proposition 1.10 in  \cite{stab}, which gives that $\G_\F$ is the underlying groupoid of $U\G_\F$. 
\end{proof}

\begin{remark}
If one wants to consider a family $\F$ that does not satisfy cancellation, one can replace $\G_\F$ by a groupoid that does satisfy cancellation (by forgetting that certain objects are isomorphic) and obtain an associated homogeneous category. We will however here for simplicity only consider families satisfying cancellation. 
\end{remark}

We end the section by showing that the homogeneous categories $U\G_\F$ considered here are not pathological in the sense that they satisfy the following {\em standardness} property: 
Let $(\C,\op,0)$ be a homogeneous category and $(A,X)$ a pair of objects in $\C$. We say that $\C$ is {\em locally standard at $(A,X)$} \cite[Definition 2.5]{stab} if 
\begin{itemize} 
\item[{\bf LS1}] The morphisms $\iota_A\op X\op \iota_X$ and $\iota_{A\op X}\op X$ are distinct in $\Hom(X,A\op X^{\op 2})$;
\item[{\bf LS2}] For all $n\ge 1$, the map $\Hom(X,A\op X^{\op n-1})\to \Hom(X,A\op X^{\op n})$ taking $f$ to $f\op \iota_X$ is injective. 
\end{itemize}

\begin{prop}\label{standard}
For any family $\F$, the category $U\G_\F$ is locally standard at any $(A,X)$.
\end{prop}

To prove this proposition, it is easiest to use an alternative description of the morphisms in the category $U\G_\F$, given by the following:

\begin{lem}\label{reform} 
The association $[X,f]\mapsto (f(X),f|_A)$ defines a one-to-one correspondence between $\Hom_{U\G_\F}(A,B)$ and the set of pairs $(H,g)$ with $H\le B$ and $g\co A\to B$ an injective homomorphism such that $B=H\x g(A)$. 
\end{lem}

\begin{proof}[Proof of Lemma~\ref{reform}]
First note that both $f(X)$ and $f|_A$ are independent of the representative of $[X,f]$, so the association is well-defined. 

Suppose that $[X,f]$ and $[Y,g]$ are morphisms from $A$ to $B$ in $U\G_\F$ satisfying that $(f(X),f|_A)=(g(Y),g|_A)$. Then $g^{-1}|_{f(X)}\circ f|_X\co X\to Y$ is an isomorphism and 
$f= g\circ ((g^{-1}|_{f(X)}\circ f|_X)\x A)$ as both maps agree on their restrictions to $X$ and $A$. Hence  $[X,f]=[Y,g]$. 

We are left to check that any $(H,g)$ is in the image. This follows from the fact that, given such an $(H,g)$, the map $H\x g\co H\x A\to B$ is an isomorphism. 
\end{proof}

\begin{proof}[Proof of Proposition~\ref{standard}]
We need to check the two axioms LS1 and LS2. 
For LS1, we need that the maps $\iota_A\x X\x \iota_X$ and $\iota_{A\x X}\x X$ from $X$ to $A\x X^2$ in $U\G_\F$ are distinct.  From the definition of the monoidal structure in $U\G_\F$ given in the proof of Proposition~1.6 of \cite{stab}, we have that $\iota_A\x X\x \iota_X=[A\x X,A\x b_{X,X}^{-1}]$ and $\iota_{A\x X}\x X=[A\x X,id_{A\x X^2}]$, where $b_{X,X}=b_{X,X}^{-1}\co X^2\to X^2$ denotes the symmetry. The fact that they are distinct then follows from the lemma as, for example, $(A\x b_{X,X}^{-1})|_{e\x e\x X}\neq id_{A\x X^2}|_{e\x e\x X}$.  

For LS2, we need to show that the map $-\x \iota_X\co  \Hom(X,A\x X^{n-1})\to \Hom(X,A\x X^{n})$ is injective. This follows again from Lemma~\ref{reform} as 
$(H,f)\x \iota_X=(H\x i_n(X),f)$ in the description of the morphisms given by the lemma, where $i_n(X)\le A\x X^n$ denotes the last $X$ factor. This association is injective.  
\end{proof}

\section{Simplicial complexes and semi-simplical sets associated to a family of groups}\label{complexes}

To a family of  groups $\F$ closed under direct
product, we associated in the previous section a category $U\G_\F$
with objects the elements of $\F$.  Using the morphism sets in this
category, the paper \cite{stab} associates to any pair of objects
$A,X\in \F$ and any $n\ge 0$,  a semisimplicial set $W_n(A,X)$ and a simplicial complex $S_n(A,X)$. In the present section, 
we recall the definitions of $S_n(A,X)$ and $W_n(A,X)$ and introduce new simplicial complexes $I_n(A,X)$ and $SI_n(A,X)$ likewise associated to $A,X\in \F$. We  
then study the relationship between these four different simplicial objects. 
To prove homological stability, we will need to show that the semisimplicial sets $W_n(A,X)$ are highly connected. This will be done in Section~\ref{connsec} in the case of the family of all right-angled Artin groups using the three simplicial complexes introduced here. We  give in the present section results that allow transfer of connectivity from one of the above spaces to another that work in a general context and that will be combined in Section~\ref{connsec} with results specific to right-angled Artin groups. 
For simplicity, we will again assume that $\F$ satisfies cancellation: 

\vspace{5mm}

\noindent
{\bf Standing assumption for the section:} $\F$ is a family of finitely generated groups, closed under direct product, and satisfying cancellation.

\bigskip

Given groups $A,X$ in $\F$, we will consider injective maps $f\co X^k\to A\x X^n$ so that there is a splitting $A\x X^n=f(X^k)\x H$ with $H$ in $\F$. As $\F$ 
satisfies cancellation, we always have that $H\cong A\x X^{n-k}$. We call such a map $f$ an {\em $\F$--split map}, and we call the pair $(f,H)$ an {\em $\F$--splitting}.

Recall that a {\em simplicial complex} $Y$ is defined from a set of vertices $Y_0$ by giving a collection of finite subsets of $Y_0$ closed under taking subsets. The subsets of cardinality $p+1$ are called the {\em $p$--simplices} of $Y$. On the other hand, a {\em semisimplicial set} $W$ is a collection of sets $W_p$ of {\em $p$--simplices} for each $p\ge 0$ related by boundary maps $d_i\co W_p\to W_{p-1}$ for each $0\le i\le p$ satisfying the simplicial identities. 
Both simplicial complexes and semisimplicial sets admit a realization, that has a copy of $\Delta^p$ for each $p$--simplex of the simplicial object. When we talk about connectivity of such objects, we always refer to the connectivity of their realization. 

\medskip

We define now three simplicial complexes and one semisimplicial set whose objects are either $\F$--split maps or $\F$--splittings.

\begin{Def}
To a pair of groups $X,A\in \F$ and a natural number $n\ge 0$, we associate the following simplicial complexes: 
\begin{description}
\item[$I_n(A,X)$] A vertex in $I_n(A,X)$ is an $\F$--split map $f\co X\to A\x X^n$. Distinct vertices $f_0,\dots,f_p$ form a $p$--simplex in $I_n(A,X)$ 
if the map $(f_0,\dots,f_p)\co X^{p+1}\to A\x X^n$ is $\F$--split. 
\item[$SI_n(A,X)$] A vertex in $SI_n(A,X)$ is an $\F$--splitting $(f,H)$ with $f\in I_n(A,X)$. Distinct vertices $(f_0,H_0),\dots,(f_p,H_p)$ form a $p$--simplex 
of $SI_n(A,X)$ if $\lgl f_0,\dots,f_p\rgl$ is a $p$--simplex of $I_n(A,X)$ and $f_i(X)\le H_j$ for each $i\neq j$. 
\item[$S_n(A,X)$] The vertices of $S_n(A,X)$ are the same as those of $SI_n(A,X)$.  Distinct vertices $(f_0,H_0),\dots,(f_p,H_p)$ form a $p$--simplex 
of $S_n(A,X)$ if there exists an $\F$--splitting $(f,H)$, with $f=(f_0,\dots,f_p)\co X^{p+1}\to A\x X^n$, such that $H_j=H\x \prod_{i\neq j}f_i(X)$ for each $j$. 
\end{description}
We moreover associate the following semisimplicial set: 
\begin{description}
\item[$W_n(A,X)$] A $p$--simplex in $W_n(A,X)$ is an $\F$--splitting $(f,H)$, with $f\co X^{p+1}\to A\x X^n$, and the $j$th face $d_j(f,H)=(f\circ d^j,H\x f(i_j))$ for 
$d^j\co X^p\to X^{p+1}$ the map skipping the $(j+1)$st factor and $i_j=\iota_{X^j}\x X\x \iota_{X^{p-j}}\co X\to X^{p+1}$. 
\end{description}
\end{Def}

Using Lemma~\ref{reform}, one checks immediately that $W_n(A,X)$ identifies with the semisimplicial set of \cite[Definition 2.1]{stab} associated to the category $U\G_\F$, and 
$S_n(A,X)$ identifies with the simplicial complex of \cite[Definition 2.8]{stab} likewise associated to $U\G_\F$.

The following proposition shows that, in the context we work with, we can always approach the connectivity of $W_n(A,X)$ via that of $S_n(A,X)$. 

\begin{prop}\label{SW} Let $\F$ be a family of groups satisfying cancellation and 
let $a,k\ge 1$.  The simplicial complex $S_n(A,X)$ is $(\frac{n-a}{k})$--connected for all $n\ge 0$ if and only if the semisimplicial set $W_n(A,X)$ is $(\frac{n-a}{k})$--connected for all $n\ge 0$. 
\end{prop}

\begin{proof}
As $U\G_\F$ is symmetric monoidal, homogeneous (Proposition~\ref{homogeneous}) and locally standard (Proposition~\ref{standard}), Proposition 2.9 of \cite{stab} yields that the semisimplicial sets $W_n(A,X)$ satisfy condition (A) in that paper (see \cite[Section 2.1]{stab}). The result then follows from \cite[Theorem 2.10]{stab}.
\end{proof}

Note that there is an inclusion of simplicial complexes $$S_n(A,X)\inc SI_n(A,X).$$
Indeed, the two complexes have the same set of vertices, and simplices of $S_n(A,X)$ satisfy the condition for being a simplex of $SI_n(A,X)$. There is also a forgetful map 
$$SI_n(A,X) \rar I_n(A,X).$$

Recall from \cite[Definition 3.2]{Hatcher-Wahl} that 
a {\em join complex\/} over a simplicial complex $X$ is a simplicial complex $Y$ together with a simplicial map $\pi\co Y\to X$  satisfying the following properties:
\begin{enumerate}[{(1)}]
\item $\pi$ is surjective;

\item $\pi$ is injective on individual simplices;

\item For each $p$-simplex $\s =\lgl x_0,\cdots,x_p\rgl$ of $X$ the subcomplex $Y(\s)$ of $Y$ consisting of all the $p$-simplices that project to $\s$ is the join $Y_{x_0}(\s)*\cdots *Y_{x_p}(\s)$ of the vertex sets $Y_{x_i}(\s)=Y(\s)\cap \pi^{-1}(x_i)$. 
\end{enumerate}
We say that $Y$ is a {\em complete join} over $X$ if $Y_{x_i}(\s)=\pi^{-1}(x_i)$ for each $\s$ and each $x_i$.

Join complexes usually arise via {\em labeling systems} (see \cite[Example 3.3]{Hatcher-Wahl}): a 
    labeling system for a simplicial complex $X$ is a
  collection of nonempty sets $L_x(\s)$ for each simplex $\s$ of $X$ and each vertex $x$ of $\s$, satisfying $L_x(\tau)\supset L_x(\s)$
  whenever $x\in\tau\subset \s$. One can think of $L_x(\s)$ as the set of labels of $x$ that are compatible with $\s$. 
We can use the labeling system $L$ to define a new simplicial complex $X^L$ having vertices the pairs $(x,l)$ with
  $x\in X$ and $l\in L_x(\lgl x\rgl)$. A collection of pairs 
  $((x_0,l_0),\cdots,(x_p,l_p))$ then forms a $p$-simplex of $X^L$ if and only if
  $\s=\lgl x_0,\cdots,x_p\rgl$ is a $p$-simplex of $X$ and $l_i\in L_{x_i}(\s)$ for each $i$. Then the natural map $\pi\co  X^L\to X$ forgetting the labels represents $X^L$ as a join complex over $X$.

\begin{prop}\label{SII}
The complex $SI_n(A,X)$ is a join complex over $I_n(A,X)$. 
\end{prop}

\begin{proof}
We check that $SI_n(A,X)$ can be constructed from $I_n(A,X)$ via a labeling system in the sense  described above. 
For each simplex $\s=\lgl f_0,\dots,f_p\rgl$ of $I_n(A,X)$ and each vertex $f_i$ in $\s$, we define the set of labels of $f_i$ compatible with $\s$ as 
$$L_{f_i}(\s):=\{H\le A\x X^n\,|\, (f_i,H)\in SI_n(A,X),\ f_j(X)\le H \ \textrm{for each}\ f_j\neq f_i\in\s\}.$$
These sets are non-empty because the fact that $\lgl f_0,\dots,f_p\rgl$ is a simplex of $I_n(A,X)$ implies that there exists an $\F$--splitting $(f,H)$ with
$$f=(f_0,\dots,f_p)\co X^{p+1}\to A\x X^n=H\x f(X^{p+1}).$$ Let $H_i=H\x \prod_{j\neq i}f_j(X)$. Then  $H_i\in L_{f_i}(\s)$. 
We clearly have that for any $f_i\in\tau\subset \s$,
$L_{f_i}(\tau)\supset L_{f_i}(\s)$, and $SI_n(A,X)=(I_n(A,X))^L$. 
\end{proof}

This will allow us to use results from \cite{Hatcher-Wahl} to obtain in good cases a connectivity bound for $SI_n(A,X)$ from one for $I_n(A,X)$. 

\medskip

We now show that, under one additional assumption, $S_n(A,X)$ and $SI_n(A,X)$ are isomorphic, in which case we will also get a connectivity result for $S_n$ from that of $SI_n$.

\begin{prop}\label{SIS}
Suppose that for any simplex $\lgl(f_0,H_0),\dots,(f_p,H_p)\rgl$ of  $SI_n(A,X)$, we have that $\bigcap_{i=0}^p H_i\in \F$. 
Then the inclusion $S_n(A,X)\inc SI_n(A,X)$ is an isomorphism. 
\end{prop}

\begin{lem}\label{splitx}
Suppose $A,B,A',B'$ are groups such that $A\x B=A'\x B'$ and $A'\le A$. Then $A=A'\x (B'\cap A)$.  
\end{lem}

\begin{proof}
Consider the inclusion $A'\x (B'\cap A)\to A$. This is an injective group homomorphism. Now every $a\in A\le A\x B=A'\x B'$ can be written as 
$a=a'b'$ with $a'\in A'$ and $b'\in B'$. But then $b'=(a')^{-1}a\in A$ and hence $a\in A'\x (B'\cap A)$ and the map is also surjective. 
\end{proof}

\begin{proof}[Proof of Proposition~\ref{SIS}]
Recall that $S_n(A,X)$ and $SI_n(A,X)$ have the same set of vertices, and that there is an inclusion   $S_n(A,X)\inc SI_n(A,X)$, that is simplices of $S_n(A,X)$ are also simplices in $SI_n(A,X)$. So we are left to check that simplices of $SI_n(A,X)$ are also always simplices in $S_n(A,X)$. 
So consider a $p$--simplex $\lgl (g_0,K_0),\dots,(g_p,K_p)\rgl$ of $SI_n(A,X)$. 
 We have  that $g=(g_0,\dots,g_p)\co X^{p+1}\to A\x X^n$ is split injective. 
 To show that these vertices form a $p$--simplex in $S_n(A,X)$, we need to 
find a complement $K\le A\x X^n$ for $g$ with $K\in \F$ satisfying 
that 
\begin{equation}\label{check}
K_j=K\x \prod_{i\neq j}g_i(X).
\end{equation}
Note that if $K$ satisfies (\ref{check}), it necessarily is a complement for $g$ as $K_j\x g_j(X)=A\x X^n$ for each $j$. 
Let $K=\bigcap_j K_j\le A\x X^{n}$. 
By the assumption, we have that $K\in\F$. We will now check that it satisfies (\ref{check}), which will finish the proof. 
By renaming the factors, it is enough to prove that (\ref{check}) holds for $j=0$. We do it by induction: we start with $K_0=K_0$. Suppose $r\ge 2$ and assume that we have proved that $K_0=\bigcap_{j=0}^{r-1}K_j\x g_1(X)\x\cdots\x g_{r-1}(X)$. 
We have that $$\bigcap_{j=0}^{r-1}K_j\x g_1(X)\x\cdots\x g_{r-1}(X) \x g_0(X) = K_0\x g_0(X)=A\x X^n=K_r\x g_r(X).$$ 
Now $g_r(X)\le K_j$ for all $j=0,\dots,r-1$. Applying the lemma we thus get that $\bigcap_{j=0}^{r-1}K_j=g_r(X)\x \big((\bigcap_{j=0}^{r-1}K_j)\bigcap K_r\big)$, which gives the induction step. 
\end{proof}

The following proposition will also be useful in the sequel:

\begin{prop}\label{transitive}
The action of 
$\aut(A\x X^n)$ on $A\x X^n$ induces an action on the complexes  $I_n(A,X)$ and $S_n(A,X)$  which is  transitive on the set of 
$p$-simplices for every $p$ in both cases. Moreover, the composed map $S_n(A,X)\to SI_n(A,X) \to I_n(A,X)$ is equivariant with respect to these actions. 
\end{prop}

\begin{proof}
The action is induced by post-composition by automorphisms on the split maps $f\co X^{p+1}\to A\x X^n$, and by evaluation on splittings $H\le A\x X^n$.   
The map $S_n(A,X)\to I_n(A,X)$ forgets the choice of splitting and is hence equivariant.  
For $S_n(A,X)$, transitivity of the action is axiom H1 in the homogeneous category $U\G_\F$ \cite[Definition 1.2]{stab}, which is satisfied by Proposition~\ref{homogeneous}. For $I_n(A,X)$, it follows from the corresponding fact for $S_n(A,X)$ and the fact that every simplex of $I_n(A,X)$ admits a lift in $S_n(A,X)$. 
\end{proof}

\section{RAAGs and their groups of automorphisms}\label{RAAGsec}

Now we consider the family $\F$ of all right-angled Artin groups, and give in this section a few properties that are particular to these groups and that will allow us to prove the connectivity result necessary for stability. 
In particular, we show that the family of RAAGs satisfies cancellation and give a description of the automorphism group of a direct product of RAAGs in terms of the automorphism groups of its factors. We start by recalling what a RAAG is. 

\medskip

Given a finite simplicial graph $\Ga$ one can associate a group
$A_\Ga$ with one generator $v$ for each vertex of $\Ga$ and a commuting relation $vw=wv$ for each edge $(v,w)$ in $\Ga$. 
Such a group $A_\Ga$ is called a graph group or more commonly {\em right-angled Artin group}. The main theorem of \cite{Droms} says that the graph describing 
such a group is unique in the sense that two such groups $A_\Ga$ and $A_{\Ga'}$ are isomorphic if and only if the graphs $\Ga$ and $\Ga'$ are isomorphic. 

\medskip

The next proposition says that RAAGs admit a ``prime decomposition'' with respect to direct product.

\begin{proposition}\label{uni} 
Any  RAAG $A_\Ga$ admits a maximal decomposition as 
$$A_\Ga=A_{\Ga_1}\x \cdots\x A_{\Ga_k}$$ with each $A_{\Ga_i}$ a RAAG, and this decomposition is unique up to isomorphism and permutation of the factors.    
\end{proposition}

\begin{proof}
From Droms' theorem \cite{Droms}, we have that the group  $A_\Ga$ splits as a direct product $A_{\Ga_1}\x A_{\Ga_2}=A_{\Ga_1*\Ga_2}$ if and only if the graph 
$\Ga$ is isomorphic to the join $\Ga_1*\Ga_2$. 
This reduces the proposition to the existence and uniqueness of the maximal decomposition of a finite graph as a join. 
Let $X$ be a finite simplicial graph. Since the graph is finite   there exists a maximal join decomposition $X=X_1*X_2*\dots *X_n$,  with $n< \infty$. Now suppose $X=X_1*X_2*\dots *X_n=Y_1*Y_2*...*Y_m$ are two distinct such decompositions. Let $k\ge 0$ be maximal such that there is a permutation $\sigma\in \Si_m$ with $X_i=Y_{\sigma(i)}$ for each $1\le i\le k$. By maximality, $X_{k+1}$ is distinct from the remaining $Y_i$'s. 
We have that $X_{k+1}$ cannot be contained in some $Y_i$ otherwise $Y_i$ would decompose as a proper join, a contradiction to the maximality of the join decomposition $Y_1*Y_2*...*Y_m$. Then $X_{k+1}$ must intersect non trivially $r$ of the remaining $Y_i$'s and so it itself must split as an $r$-join, a contradiction to the maximality of the decomposition $X_1*X_2*\dots *X_n$.  
\end{proof}

\begin{cor}
The family of RAAGs satisfies cancellation with respect to direct product. 
\end{cor}

We are here interested in the automorphism groups of RAAGs.  
The papers  \cite{Servatius, Laurence} establish that the automorphism group $\aut(A_\Ga)$ is generated by the following automorphisms:  
\begin{enumerate}[(1)]
\item (\emph{graph automorphisms}) automorphisms of the graph $\Ga$ via a permutation of its set of vertices $V$, 

\item (\emph{inversions}) for $v\in V$, a map sending $v \to v^{-1}$ and fixing all other generators, 

\item (\emph{transvections}) for $v\neq w\in V$ such that  $\link(v) \subseteq \st(w)$, a map sending $v \to vw$ and fixing all other generators,

\item (\emph{partial conjugations}) for $v\in V$ and $C$  a component of $\Ga\minus\st(v)$, the map sending $x\to vxv^{-1}$ for every vertex
$x$ of $C$ and fixing all other generators.
\end{enumerate}

\medskip

The next proposition builds on the work of Fullarton \cite{Fullarton} to  show how automorphisms of RAAGs interact with the direct product decomposition of a RAAG.

\begin{proposition}\label{Aut}
Given a RAAG $A_\Ga$ with maximal decomposition $$A_\Ga = \Z^d\x (A_{\Ga_1})^{i_1}\x \cdots\x (A_{\Ga_k})^{i_k}$$ 
with the $A_{\Ga_i}$'s  distinct and  not equal to $\Z$, 
we have 
\begin{align*}
\aut(A_\Ga) & \cong \Z^{d|\Ga'|}\rtimes\big(\GL_d(\Z)\x \aut(A_{\Ga'})\big)\\
& \cong \Z^{d|\Ga'|}\rtimes\big(\GL_d(\Z)\x  (\aut(A_{\Ga_1})\wr \Si_{i_1})\x\dots\x (\aut(A_{\Ga_k})\wr \Si_{i_k})\big)
\end{align*}
where $\Ga'=(*_{i_1}\Ga_1)*\dots*(*_{i_k}\Ga_k)$ with $|\Ga'|$ its number of vertices, and where $\Z^{d|\Ga'|}$ is generated by transvections $v\to vz$ for 
$v\in \Ga'$ and $z\in \Z^d$. 
\end{proposition}

\begin{proof}  The first isomorphism is given by Propositions 3.1 and 3.2 in \cite{Fullarton}. To get the second isomorphism, we are left to study $\aut(A_{\Ga'})$, the automorphism group of a RAAG $A_{\Ga'}$ with no $\Z$--factor.  
If $A_{\Ga'}$ is unfactorizable, there is nothing to show, so we assume that it is factorizable. We have an inclusion 
$$ (\aut(A_{\Ga_1})\wr \Si_{i_1})\x\dots\x (\aut(A_{\Ga_k})\wr \Si_{i_k})\inc \aut(A_{\Ga'})$$
so all we need to check is that every automorphism of $A_{\Ga'}$ comes from the left hand side. We do this by inspecting the generators in the classification recalled above. We see that type (2) automorphisms are internal to each factor, i.e.~are elements of some $\aut(A_{\Ga_i})$.  Type (3) can only be internal to a factor for $A_{\Ga'}$ because $A_{\Ga'}$ has no $\Z$--factors, and likewise for type (4) because $A_{\Ga'}$  is a 
  direct product. Finally, type (1) automorphisms, the graph automorphisms, satisfy 
$$\aut(\Ga')=\aut\big((*_{i_1}\Ga_1)*\dots*(*_{i_k}\Ga_k)\big)= (\aut(\Ga_1)\wr \Si_{i_1})\x\dots\x (\aut(\Ga_k)\wr \Si_{i_k}).$$
Indeed suppose $\phi$ is such a graph automorphism and let $v$ be a vertex of some copy of $\Ga_i$ and suppose that $\phi(v)$ is a vertex of a copy of some $\Ga_j$. As $\Ga_i$ is not a join, we must have that $\phi$ restricted to that $\Ga_i$ gives an injective map $\Ga_i\inc \Ga_j$. If $i=j$, this map must be an isomorphism. If not, $\Ga_j^{*i_j}$, which lies in the link of $v$, cannot be mapped to itself by $\phi$. So there must be a vertex of some $\Ga_j$ mapped to some other $\Ga_k$ with a corresponding injection induced by $\phi$.  By the pigeonhole principle, the sequence of such graph injection will end in some copy of $\Ga_i$ after finitely many steps, which then implies that in fact $\Ga_i\cong \Ga_j\cong \Ga_k\cong\dots$. Hence each $\Ga_i$ has to be mapped by such a $\phi$ to some standard copy of $\Ga_i$ in the join and the automorphism group of the join is as described. 
\end{proof}

\section{Connectivity of the simplicial complexes}\label{connsec}

In this section we show that the semisimplicial sets $W_n(A,X)$ of Section~\ref{complexes} are highly connected for any unfactorizable $X$ when $\F$ is the family of all RAAGs. 
We will treat separately the cases $X\neq\Z$ and $X=\Z$. In both cases, we will deduce this result from a computation of the connectivity of the simplicial complexes $I_n(A,X)$. 
In the first case we will show that $I_n(A,X)\cong S_n(A,X)$  while for $X=\Z$, following \cite[Section 5.3]{stab} in the case of $GL_n(R)$, we will show that $S_n(A,\Z)\cong SI_n(A,\Z)$ and use that $SI_n(A,\Z)$ is a join complex over $I_n(A,\Z)$. The connectivity of $W_n(A,X)$ will then follow using Proposition~\ref{SW}. 

The proof of connectivity of $I_n(A,X)$ when $X\neq \Z$ is a ``coloring argument'', while for $X=\Z$, 
we follow closely the work of Maazen \cite{Maazen}. 
The semisimplicial set $W_n(e,\Z)$ is essentially already in the work of Charney \cite{Charney84} under the name $SU(\Z^n)$. Charney's proof of connectivity can be adapted to the present setting and yields the same connectivity as we get.

\subsection*{Case $X\neq \Z$}

The main result of the section is the following: 

\begin{thm}\label{WnotZ}
Let $A,X$ be RAAGs such that $X\neq \Z$. Then the semisimplicial set $W_n(A,X)$ is $(n-2)$--connected. 
\end{thm}

The proof of the theorem will use the following: 

\begin{prop} \label{facts}
Let $A,X$ be RAAGs such that $X\neq \Z$ is unfactorizable and $A$ has no direct factor $X$. Then RAAG-split maps $f\co X^p\to A\x X^n$ have unique complements. Moreover, the complexes $I_n(A,X)$ and $S_n(A,X)$ are  isomorphic.
\end{prop}

\begin{proof}
The map $S_n(A,X)\to I_n(A,X)$ forgetting the chosen complements is surjective. To show that it is also injective, it is enough to check that it is injective on vertices. 
Hence the first part of the statement in the proposition in the case $p=1$  implies the second. 

By Proposition~\ref{transitive}, it is enough to check the uniqueness of complements for the standard $p$--simplex $\s_p=\lgl f_{n-p},\dots,f_{n}\rgl$ for each $p$, with $f_j\co X\to A\x X^n$ including $X$ as the $j$th $X$-factor.  The standard simplex $\s_p$ admits the subgroup $H_p=A\x X^{n-p-1}\x e\le A\x X^{n-p-1}\x X^{p+1}$ as complement. Again by Proposition~\ref{transitive}, any other complement for $\s_p$ can be obtained from $H_p$ by acting by an automorphism of $A\x X^n$  fixing the last $p+1$ factors $X$. But from the description of the automorphisms (Proposition~\ref{Aut}), we see that $H_p$ is fixed by all such automorphisms and hence $H_p$ is the only possible complement. 
\end{proof}

Let $X=\langle v_1, v_2,\ldots v_r\rangle$ where $v_1,\dots, v_r$ is the standard vertex generating set of $X$. A vertex $f \in I_n(X,
  A)$ is determined by the tuple $f=(f(v_1),\ldots, f(v_r))\in (A\x X^n)^r$. 
Write $X^n=X_1\x \cdots\x X_n$. 
 By Proposition~\ref{transitive}, we can write $f=\phi\circ f_1$ for $f_1\co X\to A\x X^n$ the inclusion as $X_1$, the first $X$-factor, and $\phi\in \aut(A\x X^n)$.  
By Proposition~\ref{Aut},  it follows that $f(v_i)=wz_{1}^{i_1}\ldots z_d^{i_d}$, for $z_{1},\dots,z_d$
  the canonical generators of $\Z^d\le A$ and $w\in X_j$ for some $j$ independent of $i$. We say that $j$ is the  \emph{color} of $f$. Note that the unique complement of $f$ is  $H= A\times  \prod_{i\neq j}(X_i)$ if the color of $f$ is $j$.

\begin{lemma}\label{Icol}
Let $A,X$ be RAAGs such that $X\neq \Z$ is unfactorizable and  is not a factor in $A$.
Then vertices $f_0,\dots,f_p\in I_n(A,X)$ form a simplex if and only if  the $f_i$'s have distinct colors in the above sense. 
\end{lemma}

\begin{proof}
Simplices of $I_n(A,X)$ have this property by Proposition~\ref{transitive} and Proposition~\ref{Aut}. Conversely, suppose $f_0,\dots,f_p$ are vertices of $I_n$ of distinct colors. Then the map $f=(f_0,\dots,f_p)\co  X^{p+1}\to A\x X^n$ 
is an injective homomorphism and  $$H=A\x \prod_{i\neq col(f_0),\dots,col(f_p)} X_i$$ is a complement for it. Hence $(f_0,\dots,f_p)$ is a $p$--simplex of $I_n$. 
\end{proof}

Recall that a simplicial complex $S$ is {\em Cohen-Macaulay} of dimension $n$ if it has dimension $n$, is $(n-1)$--connected, and the link of any $p$--simplex in $S$ is $(n-p-2)$--connected.

\begin{proposition}\label{Inconn} Let $A,X$ be RAAGs such that $X\neq \Z$ is unfactorizable and it is not a factor in $A$.
Then the  simplicial complex $I_n(X, A)$ is Cohen-Macaulay of dimension $n-1$. In particular, it is $(n-2)$-connected.
\end{proposition}

\begin{proof} Consider the map $\pi\co I_n(A,X)\to \De^{n-1}$ taking a vertex $f$  to its color. This is a simplicial map which exhibits $I_n(A,X)$ as a complete join over $\De^{n-1}$ in the sense of \cite[Definition 3.2]{Hatcher-Wahl} (see also Section~\ref{complexes}). Indeed, this map is surjective as well as injective in individual simplices. Also, for every simplex $\s=\lgl i_0,\dots,i_p\rgl$ in $\De^{n-1}$,  we have that 
$\pi^{-1}(\s)=\pi^{-1}(i_0)*\dots*\pi^{-1}(i_p)$ as vertices of $I_n$ form a simplex if and only if they have different colors by the lemma. The result is then a direct application of \cite[Proposition 3.5]{Hatcher-Wahl} and the fact that $\De^{n-1}$ is Cohen-Macaulay of dimension $n-1$. 
 \end{proof} 

We are finally ready to prove the main result of the section: 

\begin{proof}[Proof of Theorem \ref{WnotZ}] If $A\cong A'\x X^k$ for some $k>0$, we replace $W_n(A,X)$ by the isomorphic complex $W_{n+k}(A',X)$. Hence we may assume that $A$ has no $X$--factor.  
By Proposition~\ref{Inconn}, we have that $I_n(A,X)$ is $(n-2)$--connected for all $n\ge 0$. Hence by Proposition~\ref{facts}, the same holds for $S_n(A,X)$. 
Finally by Proposition~\ref{SW} with $k=1$ and $a=2$, we have that the same also holds for $W_n(A,X)$. 
\end{proof}

\subsection*{Case $X=\Z$}

The main result of the section is the following:

\begin{thm}\label{WnZ}
Let $A$ be a RAAG. Then the semisimplicial set $W_n(A,\Z)$ is $(\frac{n-3}{2})$--connected. 
\end{thm}

\begin{lem}\label{InZ}
$I_n(A,\Z)\cong I_n(e,\Z)$ for any RAAG $A$ with no $\Z$--summand. 
\end{lem}

\begin{proof}
There is an inclusion $\alpha\co  I_n(e,\Z)\to I_n(A,\Z)$ induced by composing maps to $\Z^n$ with the canonical inclusion $\Z^n\to A\x\Z^n$, given that a complement in $A\x \Z^n$ can be obtained from a complement in $\Z^n$ by crossing with $A$.
 The map $\alpha$ is simplicial 
and  injective, and we claim that it is also surjective. Indeed, by Proposition~\ref{transitive}, vertices of $I_n(\Z,A)$ are maps $f\co \Z\to A\x\Z^n$ that can be written as compositions 
$f=\phi\circ f_1$ for $f_1$ the canonical inclusion as first $\Z$--factor and $\phi$ an automorphism of $A\x \Z^n$. Now $f_1$ is in the image of  $\alpha$, and 
by the classification of the automorphisms (or Proposition~\ref{Aut}), we can see that $f=\phi\circ f_1$ still is in the image of $\alpha$: as automorphisms of $A\x \Z^n$ take $\Z^n$ to itself, the map $f$ has image in $\Z^n$. Moreover,  a complement for $f$ is of the form $A'\x H$ with $A'\cong A$ and $H\subset \Z^n$ a complement for $\al^{-1}(f)$. Likewise, if vertices form a simplex in $I_n(A,\Z)$, they will also form a simplex in $I_n(e,\Z)$. 
\end{proof}

\begin{prop}\label{InCM}
$I_n(A,\Z)$ is Cohen-Macaulay of dimension $n-1$. 
\end{prop}

Using the lemma, 
one can almost deduce the result from Corollary III.4.5 in \cite{Maazen}, though Maazen works with posets instead of simplicial complexes, and checks the vanishing of the homology groups instead of the homotopy groups. The proof adapts to our situation without any difficulty. We give it here for completeness. 

\begin{proof}
From Lemma~\ref{InZ}, we may assume that $A$ is the trivial group. For the rest of the proof, we write $I_n$ for $I_n(e,\Z)$.
We have that $I_n$ has dimension $n-1$. We need to show that it is $(n-2)$--connected, and that the link of any $p$--simplex $\s$ is $(n-p-3)$--connected for every $p\ge 0$.  Allowing $\s$ to be an empty ``$(-1)$--simplex'', we can also, and will, consider the connectivity of  $I_n$ itself  as being that of such a link. 
The link of a $p$--simplex  is non-empty whenever $n-p-2\ge 0$ and it has dimension $n-p-2$. 
We prove that the connectivity holds for each link by induction on the pair of dimensions $(\dim(\link),\dim(\s))$ in lexicographic order. The cases of $\dim(\link)=n-p-2\le 0$ are trivial as a non-empty space is $(-1)$--connected, and the empty space is $(-2)$--connected (which is defined as a non-condition). So we fix $n>p\ge -1$ with $n-p-2\ge 1$ and we assume that we have proved that $\link_{I_m}(\s)$ is $(m-k-3)$--connected for every $k$--simplex $\s$ of $I_m$ with $m>k\ge -1$ and $m-k-2\le n-p-2$, with $k<p$ if $m-k-2=n-p-2$. Let $\s$ be a $p$--simplex of $I_n$. By Proposition~\ref{transitive}, we may assume that $\s=\s_p=\lgl e_{n-p},\dots,e_n\rgl$ is the last $p+1$ standard generators in $\Z^n$, where we identify a map $f\co \Z\to \Z^n$ with the element $f(1)\in \Z^n$. We will show that $\link_{I_n}(\s)$ is $(n-p-3)$--connected, which will give the induction step and prove the result. 

\medskip

A vertex $v$ in $\link_{I_n}(\s_p)$ is given  as a $n$-tuple of integers $v=((v)_1,\dots,(v)_n)$. 
We filter the link using the absolute value of the last coordinate: let $$\OO_q:=\lgl v\in \link_{I_n}(\s_p)\ |\ |(v)_n|\le q\rgl \ \subset \ \link_{I_n}(\s_p)$$
i.e.~$\OO_q$ is the full subcomplex of $\link_{I_n}(\s_p) $ on the vertices whose last coordinate in $\Z^n$ has absolute value at most $q$. 
If $p\ge 0$, we have that $\OO_0\cong \link_{I_{n-1}}(\s_{p-1})$, the link of the last $p$ generators of $\Z^{n-1}$ in $I_{n-1}$. Indeed, if $v_1,\dots, v_k\in \Z^n$ have their last coordinate equal to 0, then 
$\lgl v_1,\dots,v_k,e_{n-p},\dots,e_n\rgl$ is a simplex of $I_n$ if and only if $\lgl \bar v_1,\dots,\bar v_k,e_{n-p},\dots,e_{n-1}\rgl$ is a simplex of $I_{n-1}$ for $\bar v_i\in \Z^{n-1}$ the first $n-1$ coordinates of $v_i$. Hence by induction, $\OO_0$ is $(n-p-3)$--connected in that case. 
If $p=-1$, $\OO_0\cong I_{n-1}$ is $(n-3)$--connected by induction. 
We will show that $\OO_1$ is $(n-2)$--connected when $p=-1$, i.e.~also $(n-p-3)$--connected. Then we will show that in both cases, for every $q\ge 0$, if $\OO_q$ is $(n-p-3)$--connected, then so is $\OO_{q+1}$. This will prove the result given that $\OO_0$ (or $\OO_1$ if $p=-1$) is $(n-p-3)$--connected, as by compactness,  any map from a sphere into the link will have image in $\OO_q$ for some $q\ge 1$. 

\medskip

We start by showing that $\OO_1$ is $(n-2)$--connected when $p=-1$. Recall that in this case $\link_{I_n}(\s_p)=I_n$. 
We can construct  $\OO_1$ from $\OO_0$ by attaching successively the vertices $v\in I_n$ with $|(v)_n|=1$, along their link in $\OO_0$, then edges formed by such vertices along their links in the newly formed complex, and so on. Explicitly this gives
$$\OO_1\  =\ \ \ \OO_0\bigcup_{v_1\in \OO_1\minus \OO_0}\hspace{-2mm} C(\lgl v_1\rgl) \hspace{-2mm}\bigcup_{\lgl v_1,v_2\rgl \subset \OO_1\minus\OO_0}\hspace{-5mm} C(\lgl v_1,v_2\rgl)\ \ \ \ \dots \hspace{-3mm}\bigcup_{\lgl v_1,\dots,v_{n}\rgl \subset \OO_1\minus \OO_0}\hspace{-7mm} C(\lgl v_1,\dots,v_{n}\rgl) $$  
where $C(\lgl v_1,\dots,v_k\rgl)=\lgl v_1,\dots,v_k\rgl * (\link_{I_n}(\lgl v_1,\dots,v_k\rgl )\cap \OO_0)$, attached successively along 
$L_k:=\del \lgl v_1,\dots,v_k\rgl * (\link_{I_n}(\lgl v_1,\dots,v_k\rgl )\cap \OO_0)$. 
For $k=1$, we have that $\link_{I_n}(\lgl v_1\rgl )\cap \OO_0=\OO_0$ as the last coordinate of $v_1$ is  $\pm 1$. Hence this link is $(n-3)$--connected. 
Now pick a vertex $v\in\OO_1\minus \OO_0$. We can write 
$$\OO_0\!\!\bigcup_{v_1\in \OO_1\minus \OO_0}\!\!\!\! C(\lgl v_1\rgl) \  =\ \ \ \st(v)\hspace{-3mm}\bigcup_{\begin{array}{c} v_1\in \OO_1\minus \OO_0\\ v_1\neq v\end{array}}\hspace{-7mm}C(\lgl v_1\rgl)$$
where $\st(v)=\OO_0*v$ is the star of $v$ within this complex.  
It follows that 
this second stage of the filtration 
is $(n-2)$-connected, being homotopic to a wedge of suspensions of $(n-3)$--connected spaces. 

For $k>1$, we have again that all the $v_i$'s have last coordinate $\pm 1$. Let $\epsilon_i=+1$ if $(v_i)_n$ and $(v_1)_n$ have the same sign, and $-1$ otherwise. Then
$$\link_{I_n}(\lgl v_1,\dots,v_k\rgl )\cap \OO_0\cong \link_{I_{n-1}}(\lgl \overline{v_2-\epsilon_2v_1},\dots,\overline{v_k-\epsilon_kv_1} \rgl ),$$ with $ \overline{v_i-\epsilon_iv_1}$ denoting as above the first $n-1$ coordinates of this vector, noting that its last coordinate is zero. Indeed, a simplex $\lgl w_1,\dots,w_q\rgl$ is in the first link if and only if the $w_j$'s have 0 as last coordinate and 
$\lgl w_1,\dots,w_q,v_1,\dots,v_k\rgl$ is a partial basis of $\Z^n$, which is the case if and only if $\lgl w_1,\dots,w_q,v_1,v_2-\eps_2v_1\dots,v_k-\eps_kv_1\rgl$ is a partial basis of $\Z^n$, which is the case if and only if $\lgl \overline{w_1},\dots,\overline{w_q},\overline{v_2-\eps_2v_1}\dots,\overline{v_k-\eps_kv_1}\rgl$ is a partial basis of $\Z^{n-1}$. 
Hence this link is $(n-k-2)$--connected by induction. So the space $L_k$, along which the cone $C(\lgl v_1,\dots,v_k\rgl)$ is attached, is $(n-2)$--connected. Hence attaching each $C(\lgl v_1,\dots,v_{k}\rgl)$ keeps the space $(n-2)$--connected. 

\medskip

We are left to show that if $\OO_q$ is $(n-p-3)$--connected then $\OO_{q+1}$ is also $(n-p-3)$--connected, where $\OO_q$ is now the $q$th filtration of $\link_{I_n}(\s_p)$ without any special assumption on $p$. We construct $\OO_{q+1}$ from $\OO_q$ by successively attaching the missing vertices, edges, and so on, just like we constructed $\OO_1$ from $\OO_0$ above: 
$$\OO_{q+1}=\OO_q\hspace{-3mm}\bigcup_{v_1\in \OO_{q+1}\minus \OO_q} \hspace{-3mm}C(\lgl v_1\rgl) \hspace{-4mm}\bigcup_{\lgl v_1,v_2\rgl \subset \OO_{q+1}\minus\OO_q} \hspace{-6mm}C(\lgl v_1,v_2\rgl)\ \ \ \ \ \dots \hspace{-6mm}\bigcup_{\lgl v_1,\dots,v_{n-p-1}\rgl \subset \OO_{q+1}\minus \OO_q}\hspace{-10mm} C(\lgl v_1,\dots,v_{n-p-1}\rgl).$$  
Again we need to compute the connectivity of the link of $\lgl v_1,\dots,v_k\rgl$ in $\link_{I_n}(\s_p)$ intersected with $\OO_q$. This link is a subcomplex of $\link_{\link_{I_n}(\s_p)}(\lgl v_1,\dots,v_k\rgl )$. By Proposition~\ref{transitive}, this last link is isomorphic to $\link_{I_n}(\s_{p+k})$, which by assumption is $(n\!-\!p\!-\!k\!-\!3)$--connected. 

Let  $\kappa\co \Z\to \Z$ be a map satisfying $\kappa(z)=0$ if  $|z|<q+1$, and  $|z-\kappa(z)(q+1)|<q+1$ for $|z|\ge q+1$. We have $(v_1)_n=\eps_1 (q+1)$ for $\eps_1=\pm 1$.  
Now define 
$$\pi\co \link_{\link_{I_n}(\s_p)}(\lgl v_1,\dots,v_k\rgl )\rar \link_{\link_{I_n}(\s_p)}(\lgl v_1,\dots,v_k\rgl )\cap \OO_q$$
to be the map taking a vertex $w$ to $w\!-\!\eps_1\kappa((w)_n)v_1$. Then $w\!-\!\eps_1\kappa((w)_n)v_1\in \OO_q$ and lies in $\link_{\link_{I_n}(\s_p)}(\lgl v_1,\dots,v_k\rgl )$ if $w$ was in that link. Moreover $\pi$ is simplicial and defines a retraction. 
It follows that $\link_{\link_{I_n}(\s_p)}(\lgl v_1,\dots,v_k\rgl )\cap \OO_q$ is also at least $(n-p-k-3)$--connected. Hence attaching $C(\lgl v_1,\dots,v_{k}\rgl)$ along 
$\del \lgl v_1,\dots,v_k\rgl * (\link_{\link(\s_p)}(\lgl v_1,\dots,v_k\rgl )\cap \OO_q)$ does not change the connectivity as the latter space is at least $(n-p-3)$--connected. The result follows.
\end{proof}

\begin{remark}
The existence of the function $\kappa$ used in the proof makes $\Z$, together with the absolute value, a {\em Euclidean ring}. Maazen's proof of the above statement was written, and works in the more general  context of Euclidean rings.  
\end{remark}

\begin{prop}\label{intersect}
Let $A$ be a RAAG with no $\Z$--summand and let $\lgl (K_0,f_0),\dots,(K_p,f_p)\rgl$ be a simplex of $SI_n(A,\Z)$. Then $\bigcap_j K_j\cong A\x \Z^m$ for some $m$. 
\end{prop}

\begin{proof} Suppose that $A$ has generators $a_1,\dots,a_r$ and $\Z^n$ has generators $z_1,\dots,z_n$. 
We know that each $K_i$ can be obtained from the standard $A\x\Z^{n-1}\x e\le A\x\Z^n$ by applying an automorphism. From the description of the automorphisms of $A\x\Z^n$, it follows that $K_i$ is generated by $a_1w_1,\dots,a_rw_r,t_1,\dots,t_{n-1}$ for some $w_j,t_j\in \Z^n$.  Moreover, we know that $\Z^n$ is generated by $t_1,\dots, t_{n-1},f_i$, where $f_i:=f_i(1)$. Hence we can rewrite the generators of $K_i$ as $a_1\,f_i^{\,m_{i,1}},\dots,a_r\,f_i^{\,m_{i,r}},t_1,\dots,t_{n-1}$ for some $m_{i,1},\dots,m_{i,r}\in \Z$. 

As $f_i\in K_j$ whenever $i\neq j$ and $a_k\,f_i^{\,m_{i,k}}\in K_i$, we have that each $a_k\,f_0^{\,m_{0,k}}f_1^{\,m_{1,k}}\dots f_p^{\,m_{p,k}}$ lies in $\bigcap_i K_i$. Let $A'\cong A$ denote the subgroup of $A\x \Z^n$ generated by the elements $a_k\,f_0^{\,m_{0,k}}f_1^{\,m_{1,k}}\dots f_p^{\,m_{p,k}}$'s. We have that $A\x \Z^n=A'\x \Z^n$. 
We want to show that $\bigcap_i K_i=A'\x (\bigcap_i (K_i\cap \Z^n))$. The right side is included in the left side, so all we need to show is that the left side is included in the right side. Let $x\in\bigcap_i K_i$ be some element. As $x\in K_i$, we can write it as $x=x'_ix''_i$ with $x'_i\in A'$ and $x_i''\in K_i\cap \Z^n$. Now these different expressions of $x$ are all equal, and all live in $A\x \Z^n=A'\x \Z^n$. It follows that $x'_i=x'_j$ and hence $x''_i=x''_j$ for each $i,j$. It follows that $x''_0\in \bigcap_i(K_i\cap \Z^n)$  and $x=x_0'x_0''\in A'\x (\bigcap_i (K_i\cap \Z^n))$. As $\bigcap_i (K_i\cap \Z^n)\le \Z^n$, the result follows. 
\end{proof}

\begin{proof}[Proof of Theorem~\ref{WnZ}]
Just as in the proof of Theorem~\ref{WnotZ}, we may assume that $A$ has no $\Z$--factor. 
From Proposition~\ref{SII}, we have that $SI_n(A,\Z)$ is a join complex over $I_n(A,\Z)$ (in the sense of \cite[Definition 3.2]{Hatcher-Wahl}, see also Section~\ref{complexes}).
As $I_n(A,\Z)$ is Cohen-Macaulay of dimension $n-1$ by Proposition~\ref{InCM}, Theorem 3.6  of \cite{Hatcher-Wahl} gives that $SI_n(A,\Z)$ is $(\frac{n-3}{2})$--connected for all $n\ge 0$. By Proposition~\ref{intersect}, the hypothesis of  Proposition~\ref{SIS} is satisfied, and hence $S_n(A,\Z)$ is isomorphic to $SI_n(A,\Z)$. So the connectivity also holds for $S_n(A,\Z)$. Hence by  Proposition~\ref{SW} with $k=2$ and $a=3$, we have that the same also holds for $W_n(A,X)$. 
\end{proof}

\section{Stability theorem}\label{stabsec}

We consider in this section the family of all right-angled Artin groups with $\C=U\G_{Raags}$ the associated homogeneous category, as defined in Section~\ref{sec1}. 
Let $A,X$ be RAAGs, and denote by $\C_{A,X}$ the full subcategory of $\C$ on the objects $A\x X^n$ for all $n\ge 0$. 
Recall from \cite[Section 4.2]{stab} the {\em lower suspension} functor $\Si_X\co \C_{A,X}\to \C_{A,X}$ taking $A\x X^n$ to $A\x X^{n+1}$ and a morphism $f\co A\x X^n\to A\x X^k$ to the composition $(b_{X,A}\x X^k)\circ (X\x f) \circ (b_{X,A}^{-1}\x X^n)$, where $b_{X,A}\co X\x A \to A\x X$ denotes the symmetry. 
Recall from \cite[Definition 4.10]{stab} that a functor $$F\co \C_{A,X}\rar \Z\operatorname{-Mod}$$
is a {\em coefficient system of degree $r$ at $N$ } if the kernel of the suspension map $F\to F\circ \Si_X$ is trivial when evaluated at $A\x X^{n}$ with $n\ge N$, and the cokernel is of degree $r-1$ at $N-1$, with degree $-1$ at $N$ meaning taking the value 0 at $A\x X^{n}$ whenever $n\ge N$. 
In particular, constant coefficient systems are of degree 0 at 0. 
A coefficient system $F$ is {\em split} if it splits as a functor.  

Applying the main results of \cite{stab} to our situation, we get the following stability theorem:

\begin{thm}\label{allbounds}
Let $A,X$ be RAAGs with $X$ unfactorizable and $F\co \C_{A,X}\to \Z\operatorname{-Mod}$ a coefficient system of degree $r$ at $N$. Let $n>N$. 
Then the map $\aut(A\x X^n)\to \aut(A\x X^{n+1})$ taking an automorphism $f$ to $f\x X$  induces a map 
$$H_i(\aut(A\x X^n);F(A\x X^n))\rar H_i(\aut(A\x X^{n+1});F(A\x X^{n+1}))$$ 
which is surjective for all $i\le \frac{n-1}{2}-r$ and an isomorphism for all  $i\le\frac{n-3}{2}-r$. If the coefficient system is split, this range improves to 
 $i\le \frac{n-r-1}{2}$ for surjectivity and  $i\le \frac{n-r-3}{2}$ for injectivity, and if the coefficient system is constant, the isomorphism holds for  $i\le \frac{n-2}{2}$.

Moreover, let $\aut'(A\x X^n)$ denote the commutator subgroup of $\aut(A\x X^n)$. Let $n>2N$. Then the map 
$$H_i(\aut'(A\x X^n);F(A\x X^n))\rar H_i(\aut'(A\x X^{n+1});F(A\x X^{n+1}))$$ 
is surjective for all $i\le \frac{n-2}{3}-r$ and an isomorphism for all  $i\le\frac{n-5}{3}-r$. If the coefficient system is split, this range improves to 
 $i\le \frac{n-2r-2}{3}$ for surjectivity and  $i\le \frac{n-2r-5}{3}$ for injectivity, and if the coefficient system is constant, the isomorphism holds for  $i\le \frac{n-4}{3}$. 

If $X\neq \Z$--summand, one can replace $n$ by $n+1$ in all the bounds of the theorem.  
\end{thm}

\begin{proof}
The category $\C=U\G_{Raags}$ is symmetric monoidal homogeneous by Proposition \ref{homogeneous}, and hence pre-braided and locally homogeneous at any $(A,X)$ in the sense of \cite[Definition 1.1,1.4]{stab}. 
Theorem \ref{WnotZ} gives that it satisfies LH3 with slope 2 at $(A,X)$ for all $A$ and all irreducible $X\neq \Z$ (\cite[Definition 2.2]{stab}), and Theorem \ref{WnZ} that  it satisfies LH3 with slope 2 at $(A\x \Z ,\Z)$ for all $A$. 
The result then follows from Theorems 3.1, 3.4 and 4.20 of \cite{stab} for $(A,X)$ with $X$  unfactorizable not equal to $\Z$, and for $(A\x \Z,\Z)$, using the argument of \cite[Corollary 3.9]{stab} for the second part of the statement. 
\end{proof}

Theorems~\ref{maincst} and \ref{mainab} are obtained by applying the above theorem to constant coefficient systems for each irreducible factor of  $B=X_1\x \cdots\x X_k$. 
Theorem~\ref{mainpol} is obtained likewise applying the theorem to the coefficient system defined by the abelianization, noting that this is a split coefficient system of degree 1 at 0.

Note that van der Kallen obtains better bounds for $\GL_n(\Z)$, which is the case when $A$ is the trivial group and $B=X=\Z$. Most particularly, for Theorem~\ref{mainab}, his bound has slope 2 instead of slope 3 as we have (see \cite[Theorem 4.6]{vdk} or Proposition~\ref{abeliangp}). However his argument does not obviously extend to all RAAGs. (The argument of van der Kallen is explained at the end of Section 5.3 in \cite{stab}.)

\subsection{Finitely generated abelian groups}

Homological stability for the automorphism groups of finitely generated abelian groups under taking direct product can be deduced directly from existing results in the literature, without needing to prove new connectivity results. We give the exact statement and its proof here for completeness.

\begin{prop}\label{abeliangp}
Let $G$ be a finitely generated abelian group. The homomorphism $\aut(G^n)\to \aut(G^{n+1})$ taking an automorphism $f$ of $G^n$ to the automorphism $f\x G$ of $G^{n+1}$ fixing the last factor, induces maps 
$$H_i(\aut(G^n);\Z)\rar  H_i(\aut(G^{n+1});\Z)$$
and 
$$H_i(\aut'(G^n);\Z)\rar  H_i(\aut'(G^{n+1});\Z)$$
which are surjective  for all $i\le \frac{n}{2}$ and isomorphisms for $i\le \frac{n-1}{2}$.
\end{prop}

(As in the introduction,  $\aut'(G^n)$ denotes the commutator subgroup of $\aut(G^n)$.) 

\begin{proof}
A finitely generated abelian group $G$ is a $\Z$--module, and the automorphism group $\aut(G^n)$ is isomorphic to $\GL_n(\operatorname{End}(G))$, for  $\operatorname{End}(G)$ its ring of endomorphisms.  
By  \cite[Theorem 3.4]{War80}, $\operatorname{End}(G)$ has 2 in its stable range (in the terminology of \cite{War80}, see Definition~1.5 in that paper), i.e.~satisfies Bass' condition $SR_3$ or has $sdim=1$ in the terminology of \cite[Section 2.2]{vdk}. The result then follows from  \cite[Theorem 4.11]{vdk} using the fact that $\GL_n'(R)$ is isomorphic to its subgroup of elementary matrices (Whitehead's lemma). 
\end{proof}

Theorem 5.6 of \cite{vdk} and Theorem 5.10 of \cite{stab} can likewise be applied to show that homological stability for the groups $\aut(G^n)$ with $G$ finitely generated abelian also holds with polynomial twisted coefficients.

\bibliographystyle{plain}
\bibliography{math}

\begin{thebibliography}{10}

\bibitem{Charney84}
Ruth Charney.
\newblock On the problem of homology stability for congruence subgroups.
\newblock {\em Comm. Algebra}, 12(17-18):2081--2123, 1984.

\bibitem{gln}
Ruth~M. Charney.
\newblock Homology stability for {${\rm GL}_{n}$} of a {D}edekind domain.
\newblock {\em Invent. Math.}, 56(1):1--17, 1980.

\bibitem{djament}
Ga{\"e}l Collinet, Aur{\'e}lien Djament, and James~T. Griffin.
\newblock Stabilit\'e homologique pour les groupes d'automorphismes des
  produits libres.
\newblock {\em Int. Math. Res. Not. IMRN}, (19):4451--4476, 2013.

\bibitem{Droms}
Carl Droms.
\newblock Isomorphisms of graph groups.
\newblock {\em Proc. Amer. Math. Soc.}, 100(3):407--408, 1987.

\bibitem{Fullarton}
N.~J. {Fullarton}.
\newblock {On the number of outer automorphisms of the automorphism group of a
  right-angled Artin group}.
\newblock {\em Math. Res. Lett. (to appear)}.

\bibitem{Gra76}
Daniel Grayson.
\newblock Higher algebraic {$K$}-theory. {II} (after {D}aniel {Q}uillen).
\newblock In {\em Algebraic {$K$}-theory ({P}roc. {C}onf., {N}orthwestern
  {U}niv., {E}vanston, {I}ll., 1976)}, pages 217--240. Lecture Notes in Math.,
  Vol. 551. Springer, Berlin, 1976.

\bibitem{hatcher}
Allen Hatcher.
\newblock Homological stability for automorphism groups of free groups.
\newblock {\em Comment. Math. Helv.}, 70(1):39--62, 1995.

\bibitem{HatVog}
Allen Hatcher and Karen Vogtmann.
\newblock Cerf theory for graphs.
\newblock {\em J. London Math. Soc. (2)}, 58(3):633--655, 1998.

\bibitem{Hatcher-Wahl}
Allen Hatcher and Nathalie Wahl.
\newblock Stabilization for mapping class groups of 3-manifolds.
\newblock {\em Duke Math. J.}, 155(2):205--269, 2010.

\bibitem{Hirshon0}
R.~Hirshon.
\newblock The cancellation of an infinite cyclic group in direct products.
\newblock {\em Arch. Math. (Basel)}, 26:134--138, 1975.

\bibitem{Hirshon1}
R.~Hirshon.
\newblock Some cancellation theorems with applications to nilpotent groups.
\newblock {\em J. Austral. Math. Soc. Ser. A}, 23(2):147--165, 1977.

\bibitem{john}
Francis E.~A. Johnson.
\newblock Automorphisms of direct products of groups and their geometric
  realisations.
\newblock {\em Math. Ann.}, 263(3):343--364, 1983.

\bibitem{Laurence}
Michael~R. Laurence.
\newblock A generating set for the automorphism group of a graph group.
\newblock {\em J. London Math. Soc. (2)}, 52(2):318--334, 1995.

\bibitem{Maazen}
Henrik Maazen.
\newblock {Homology stability for the general linear group}.
\newblock www.staff.science.uu.nl/$\sim$kalle101/maazen1979.pdf, 1979.

\bibitem{stab}
Oscar Randal-Williams and Nathalie Wahl.
\newblock Homological stability for automorphism groups.
\newblock arXiv:1409.3541, 2014.

\bibitem{Servatius}
Herman Servatius.
\newblock Automorphisms of graph groups.
\newblock {\em J. Algebra}, 126(1):34--60, 1989.

\bibitem{vdk}
Wilberd van~der Kallen.
\newblock Homology stability for linear groups.
\newblock {\em Invent. Math.}, 60(3):269--295, 1980.

\bibitem{War80}
R.~B. Warfield, Jr.
\newblock Cancellation of modules and groups and stable range of endomorphism
  rings.
\newblock {\em Pacific J. Math.}, 91(2):457--485, 1980.

\end{thebibliography}

\end{document}